\documentclass{article}
\usepackage[utf8]{inputenc}
\usepackage{lipsum}
\usepackage{url}
\usepackage[T1]{fontenc}
\usepackage[utf8]{inputenc}
\usepackage[english]{babel}
\usepackage{times}
\usepackage[left=1.8cm, right=1.8cm, top=3cm]{geometry}
\usepackage{amsmath}
\usepackage{amssymb}
\usepackage{amsthm}
\usepackage{hyperref}
\usepackage{graphicx}
\usepackage{subcaption}
\usepackage{caption}
\usepackage{xcolor}
\usepackage{longtable}
\usepackage{verbatimbox}

\usepackage[title]{appendix}
\usepackage{bm}
\usepackage{empheq}
\usepackage{placeins}

\usepackage{xcolor}
\hypersetup{
colorlinks,
linkcolor={red!50!black},
citecolor={blue!50!black},
urlcolor={blue!80!black}
}

\newcommand\scalemath[2]{\scalebox{#1}{\mbox{\ensuremath{\displaystyle #2}}}}

\usepackage{amsmath}
\usepackage{amssymb}
\usepackage{amsthm}
\usepackage{hyperref}
\usepackage{graphicx}
\usepackage{subcaption}
\usepackage{caption}
\usepackage{xcolor}

\usepackage{longtable}
\usepackage{verbatimbox}
\usepackage[title]{appendix}
\usepackage{bm}
\usepackage{empheq}
\usepackage{placeins}
\theoremstyle{plain}
\newtheorem{theorem}{Theorem}[section]

\newtheorem{remark}[theorem]{Remark}
\newtheorem{definition}[theorem]{Definition}
\theoremstyle{definition}
\theoremstyle{remark}
\numberwithin{equation}{section}
\newenvironment{oss}{\begin{remark} \begin{rm}}{\end{rm} \end{remark}}
\newcommand\norm[1]{\left\lVert#1\right\rVert}
\newcommand{\ti}{{\vect{t}}}
\newcommand{\di}{{\vect{d}}}
\newcommand{\rr}{{\vect{r}}}
\newcommand{\zz}{{\vect{z}}}

\theoremstyle{plain} 

\usepackage{mathtools}
\usepackage{pgfplots}
\pgfplotsset{/pgf/number format/use comma,compat=newest}

\renewcommand\epsilon{\varepsilon}
\newcommand{\R}{\mathbb{R}}

\newcommand{\vect}[1]{\boldsymbol{#1}}

\usepackage{subcaption}
\usepackage[small, sc]{titlesec}

\usepackage{xcolor}
\hypersetup{
colorlinks,
linkcolor={red!50!black},
citecolor={blue!50!black},
urlcolor={blue!80!black}
}

\usepackage[
backend=biber,
style=numeric,
sorting=nty,
maxbibnames=50,
giveninits=true
]{biblatex}
\usepackage{csquotes}
\renewbibmacro{in:}{%
 \ifentrytype{article}{}{}}
\DeclareFieldFormat[article]{volume}{\textbf{#1}}
\DeclareFieldFormat[article]{title}{\textit{#1}}
\DeclareFieldFormat[article]{journaltitle}{#1}
\DeclareFieldFormat[article]{number}{\textbf{#1}}
\renewbibmacro*{volume+number+eid}{%
  \printfield{volume}%
\setunit*{\textbf{\addcolon}}
  \printfield{number}\setunit{\addcomma\space}%
  \printfield{eid}}

\graphicspath{{figure/}}

\begin{document}

\title{\textsc{
Effects of surface tension and elasticity on critical points of the Kirchhoff-Plateau problem
}}

\author{
\textsc{Giulia Bevilacqua}\\[7pt]
\small Dipartimento di Matematica\\ 
\small Università di Pisa\\ 
\small Largo Bruno Pontecorvo 5, I–56127 Pisa, Italy\\
\small \href{mailto:giulia.bevilacqua@dm.unipi.it}{giulia.bevilacqua@dm.unipi.it}
\and
\textsc{Chiara Lonati}\\[7pt]
\small  Dipartimento di Matematica ``F. Casorati''\\
\small Università degli Studi di Pavia
\\ 
\small via Ferrata 5, I-27100 Pavia, Italy
\\
\small \href{mailto:chiara.lonati01@universitadipavia.it}{chiara.lonati01@universitadipavia.it}
}

\date{}

\maketitle

\begin{abstract}
\noindent
We introduce a modified Kirchhoff-Plateau problem adding an energy term to penalize shape modifications of the cross-sections appended to the elastic midline. In a specific setting, we characterize quantitatively some properties of minimizers. Indeed, choosing three different geometrical shapes for the cross-section, we derive Euler-Lagrange equations for a planar version of the Kirchhoff-Plateau problem. We show that in the physical range of the parameters, there exists a unique critical point satisfying the imposed constraints. Finally, we analyze the effects of the surface tension on the shape of the cross-sections at the equilibrium.
\end{abstract}

\bigskip



\bigskip

\maketitle

\section{Introduction}
\label{sec:intro}
The minimization of the area functional dates back to Plateau \cite{plateau1873experimental} who named the famous {\em Plateau problem}, which in its simplest form asks if it exists a minimal surface spanning a given curve. For the first rigorous mathematical result, we need to wait Douglas and Rad\'o \cite{rado1930problem, douglas1931solution}, who, separately, around 1930s substituted the area functional with the Dirichlet one using conformal maps. \textcolor{black}{Later on}, many mathematicians generalized such a definition for the surface, like finite perimeter sets \cite{maggi2012sets}, currents \cite{federer1960normal} and Almgren minimal sets  \cite{almgren1976existence, taylor1976structure} (we refer to the complete survey \cite{david2014should} and references therein).

A recent variant of the Plateau problem considers a $3$D object as the assigned boundary, namely a tubular neighbourhood of a given curve. The most delicate aspect is prescribing the intersection set between the minimal surface and the external boundary of the rigid body. A weak topological notion of {\em spanning} introduced by Harrison and Pugh \cite{harrison2016existence} allows for the treatment of the free-boundary problem. Indeed, using the previous definition, De Lellis, Ghiraldin and Maggi proved the existence of an optimal soap film surface \cite{DeLellis_2014} with the right regularity in the physical dimension \cite{almgren1968existence, taylor1976structure}. We also mention a recent capillarity model for soap film\textcolor{black}{s} as regions containing small volume with homotopic spanning condition \cite{maggi2019soap}.

In contrast to the Plateau problem, the {\em Kirchhoff--Plateau} problem concerns the equilibrium shapes of a system composed by a closed Kirchhoff rod spanned by an area-minimizing surface. The main difference is the presence of a competition between normal forces of the soap film and the elasticity of the rod \cite{roman2010elasto}. This requires a suitable compactness theorem since boundaries change along minimizing sequences \cite{Lussardi_2017}. In such a situation, proving the existence of minima of the energy functional under physical constraints is quite well-understood and we mention several results in this direction \cite{biria2014buckling, biria2015theoretical, bevilacqua2018soap, bevilacqua2019soap, bevilacqua2020dimensional}.

The study of critical points for the Kirchhoff-Plateau problem is still a matter of debate. Indeed, such a derivation seems to be hard for two reasons. On one hand, \textcolor{black}{since} the rod is highly constrained, \textcolor{black}{computations of suitable variations are quite hard}. On the other hand, the regularity of the contact set between rod and surface is not yet known \cite{DeLellis_2014}. Recently, in \cite{king2022plateau} the authors derived in a weak sense Euler-Lagrange equations for {\em generalized minimizers} containing a fixed volume and spanning an assigned tubular neighbourhood. For the elastic counterpart, few results are known: we  mention \cite{giusteri2016instability} for a formal derivation of the first-necessary conditions and \cite{bevilacqua2022variational, BBLM22} where the elastic rod is replaced by an elastic curve.

For all these reasons, to explicitly derive Euler-Lagrange equations for an elastic rod spanned by a soap film, we need to make some simplifications. Indeed, for a complete result, we expect to deal with Plateau's type singularities, suitable variations for the rod and the use of a Geometric Measure Theory framework: this will be the content of a forthcoming paper.

First of all, we assume higher regularity, hence the midline $\vect{r} \in \mathcal{C}^{\infty}([0,L];\R^3)$, where $L>0$ is the length of the curve. Second, we restrict our attention to a planar problem, constraining the midline to remain in the \textcolor{black}{horizontal plane $z = 0$}, so that we describe the curve using only a geometrical parameter. Then, to quantify the competition between normal forces of the soap film and the elasticity of the rod, we assume that the trace of the surface on the rod is a curve $\vect{\tilde{r}}$, \textcolor{black}{called {\em contact curve}}. Such a curve is a ``similar midline'', \textcolor{black}{indeed it is regular, $\vect{\tilde{r}} \in C^{\infty}([0,L];\R^3)$ and it lies in the same plane as the one of the midline, hence $\vect{\tilde{r}} \in \{z = 0\}$.}

In the following, we will consider energy functionals of the form
$$
\mathcal{F}[\vect{r},\Phi]:= \int_0^L f(\kappa, s)\, ds + \int_0^L \Phi_{\vect{p}}(s)\, ds + 2 \sigma \mathcal{H}^2(S_{\vect{\tilde{r}}}),
$$
where $\kappa$ is the curvature, the only geometric invariant we need to describe a planar curve, $\sigma$ is the surface tension and $S_{\vect{\tilde{r}}}$ refers to the area of the closed surface spanning the scaled midline ${\vect{\tilde{r}}}$. We define $\Phi_{\vect{p}}$ as a surface density to model \textcolor{black}{shape modifications} 
of the cross-sections appended to the midline $\vect{r}$ \cite{antman} due to normal forces of the minimal surface on the boundary of the rod \cite{roman2010elasto}. The subscript $\vect{p}$ refers to the influence of such a function to the reconstruction process of the rod and its position in $\R^3$. Moreover, choosing different shapes of the cross-sections, we aim to characterize the effect of the surface tension of the film on equilibrium configurations of the functional $\mathcal{F}$.

The paper is organized as follows. We present the general setting for the Kirchhoff-Plateau problem with variable cross-sections in Section \ref{sec:generale}.
Then, in Section \ref{sec:planar}, we introduce the necessary simplifications in the planar context to derive the first-order necessary conditions. In Section \ref{elsec}, we consider an elliptical shape for the cross-section; in Section \ref{sec:homotetia} a dilated ellipse along the horizontal axis and finally in \ref{sec:oval} an oval shape. In all the studied cases, we derive conditions on the parameters to obtain equilibrium configurations and we characterize how the surface tension modifies the shape of the cross-section.
We conclude with some remarks in Section \ref{sec:conclusions}.

\section{Setting of the problem and general energy functional}
\label{sec:generale}
In this Section, we briefly introduce the Kirchhoff-Plateau problem and an energy functional which takes into account an elastic term for the variation of the cross-section. 

\subsection{The elastic rod}

The Kirchhoff-Plateau problem asks if there exists a surface with minimal area spanning an elastic thick boundary. Following \cite{antman}, an elastic rod is a three-dimensional object described through a curve $\vect{r}$ \textcolor{black}{which we call} {\em midline} and a family of material frames appended to it. 
We assume that the midline  is inextensible, i.e. it can be parametrized through the arc-length parameter $s \in [0,L]$ where $L>0$ is the length of the curve, and that the rod is unshearable, i.e. at each point on the curve, the plane of the cross-section \textcolor{black}{$\mathcal{A}(s)\subset \mathbb{R}^2$} is orthogonal to the midline. \textcolor{black}{The material cross-section $\mathcal{A}(s)$ is, for every $s$, a compact simply connected set, such that the origin $(0,0) \in \R^2$ belongs to $\operatorname{int}(\mathcal{A}(s))$}. Under such hypotheses, the shape of the rod is uniquely determined assigning a function $s\in [0,L]\mapsto \mathcal{A}(s)$ to define the shape of the cross-sections, and the material coefficients $\kappa_1, \kappa_2,\omega$, named  {\em flexural densities} and {\em twist density} respectively, to describe mechanical properties of the rod. Now, choosing the shape of $\mathcal{A}(s)$, setting ${\rm w}:=(\kappa_1(s),\kappa_2(s), \omega(s)) \in L^2(0,L)\times L^2(0,L)\times L^2(0,L)$, the midline $\rr$ and the director field $\di$ orthogonal to the tangent vector $\vect{t}(s) = \vect{r}'(s)$ can be uniquely determined through the vectorial system of ODEs  \cite{hartman}
\begin{equation}
\label{eq:system}
	\left\{
	\begin{aligned}
	&\ti'(s) = \kappa_1(s) \di(s) + \kappa_2(s) \zz(s),\\
	&\di'(s) = -\kappa_1(s) \ti(s) + \omega(s) \zz(s),\\
	&\zz'(s) = -\kappa_2(s) \ti(s) -\omega(s) \di(s),
	\end{aligned}
	\right.
\end{equation}
equipped with the clamping conditions
$$
\begin{aligned}
&&\ti(0) = \ti_0 &&\di(0) = \di_0 &&\zz(0) = \zz_0,
\end{aligned}
$$
where $\{\ti_0(s),\di_0(s), \zz_0(s)\}$ is an orthonormal basis. By classical results on vectorial systems of ODEs \cite{hartman}, we get as solution of \eqref{eq:system} the midline $\rr(s)$ and the orthonormal frame $\{\ti(s),\di(s), \zz(s)\}$, where $\zz:= \ti \times \di$, implying that the last equation of \eqref{eq:system} is always true providing $\vect{t}$ and $\vect{d}$. Moreover, the following regularity holds
\begin{equation}
\label{eq:soluzione}
\begin{aligned}
&\rr(s) := \vect{r}_0 + \int_0^s \vect{t}(\tau)\, d\tau\; \in W^{2,2}((0,L);\R^3),&&&(\ti,\di) \in W^{1,2}((0,L);\R^3)\times W^{1,2}((0,L);\R^3),
\end{aligned}
\end{equation}
where $\vect{r}_0$ is the assigned condition to reconstruct the curve in $\R^3$ and obviously $\vect{z} = \vect{t} \times \vect{d} \in W^{1,2}((0,L);\R^3)$.
Then, the position of the points of the rod is defined through a function $\vect{p}:\Omega\to \R^3$, where
\begin{equation}
\label{eq:defp}
    \Omega := \{(s,\zeta_1,\zeta_2): s\in [0,L] \hbox{ and } (\zeta_1,\zeta_2) \in \mathcal{A}(s)\} \quad \hbox{ and } \quad \vect{p}(s, \zeta_1, \zeta_2):= \vect{r}(s) + \zeta_1 \vect{d}(s) + \zeta_2 \vect{z}(s). 
\end{equation}

Before introducing the energy, we list the constraints imposed to the rod (we refer to \cite{Schuricht_2002} for details).
\begin{enumerate}
    \item  The rod and the tangent vector are closed, hence
\begin{equation}
    \label{eq:chiusura}
    \begin{aligned}
        &\vect{r}(L)=\vect{r}(0)=\vect{r}_0, &&&\vect{t}(L)=\vect{t}(0)=\vect{t}_0.
    \end{aligned}
\end{equation}
   \item To glue the rod, we fix $\vect{d}(L) = \vect{d}_0$ and we prescribe how many times the ends of the rod can be twisted before being glued together fixing the linking number $L_0$ of the midline with a near suitable closed curve.
   \item We encode the knot type of the midline fixing a continuous mapping $\vect{\eta}: [0,L]\to \R^3$ such that $\vect{\eta}(0) = \vect{\eta}(L)$ and we require the midline to be isotopic to $\vect{\eta}$.
    \item To prevent self-intersections, we require the Ciarlet-Ne\v cas condition \cite{Ciarlet_1987}, namely
    \[
\int_{\Omega}\det \left(D\vect{p}[{\rm w}]\right)\,ds\,d\zeta_1\,d\zeta_2\leq \mathcal{L}^3(\vect{p}[{\rm w}](\Omega)).
\]
\end{enumerate}
We denote by $U$ the set of all constraints as
\begin{equation}
\label{eq:set_constraint}
U  = \left\{{\rm w}\in L^2((0,L);\R)^3: \hbox{1 - 4 hold true}\right\},
\end{equation}
which turns out to be weakly closed in $L^2((0,L);\R)^3$, see \cite{Schuricht_2002}.

Finally, we can introduce the elastic and potential energy stored in the rod. Let $f : [0,L]\times \R^3 \to \R \cup \{+ \infty\}$ be bounded from below and such that $\forall \, s \in [0,L], \,f(s,\cdot)$ is continuous and convex and $\forall\, a \in \R^3, \,f(\cdot, a)$ is measurable and coercive. Then, the elastic contribution of the rod $E_{\rm el}[{\rm w}]:U \to \R \cup \{+ \infty\}$ reads
$$
E_{\rm el}[{\rm w}]:=  \int_0^L f(s, {\rm w})\, ds.
$$
As for the potential energy of the weight, it is given by
$$
E_{\rm g}[{\rm w}]:= -\int_\Omega \rho(s, \zeta_1, \zeta_2)\vect{g} \cdot \vect{p}[{\rm w}](s, \zeta_1, \zeta_2)\, dsd\zeta_1d\zeta_2,
$$
where $\rho$ is the given mass density and $\vect{g}$ is the gravitational acceleration.
 To take into account a possible modification of the cross-section of the rod, we introduce an \textcolor{black}{additional contribution, that corresponds to shape modifications of the cross-section. Let $\Phi_{\vect{p}}:[0,L]\to [0, +\infty)$ be an energy density} such that $\Phi_{\vect{p}}$ is a continuous function; such a contribution reads
\begin{equation}
    \label{eq:E_sh_tot}
    E_{\rm sh}[\Phi]:= \int_0^L \Phi_{\vect{p}}(s)\, ds.
\end{equation}                                             
The function $\Phi_{\vect{p}}$ represents a surface density encoding an \textcolor{black}{energy contribution of the cross-section $\mathcal{A}(s)$ for all $s \in [0,L]$, that penalizes modifications of its shape}. Roughly speaking, $\Phi_{\vect{p}}$ is the integral over the cross-section $\mathcal{A}(s)$ of a suitable energy function which specifies different geometrical shapes of the cross-section. Precisely, $\Phi_{\vect{p}}$ encodes the energy density of change of shapes of the cross-sections concentrated on the midline. \textcolor{black}{The subscript $\vect{p}$ refers to the influence of such an energy on the set $\Omega$ and on modifications of the function $\vect{p}$ in \eqref{eq:defp}.  Hence, to define the position of the rod in $\R^3$ through $\vect{p}(\Omega)$, we need to take care of the shape of $\mathcal{A}(s)$ (remember that \eqref{eq:defp} contains $\zeta_1, \zeta_2 \in \mathcal{A}(s)$): adding a shape modification term to the functional implies that the configuration map $\vect{p}$, the elastic energy and the gravitational energy of the rod are influenced by $\Phi_{\vect{p}}$.}
Then, the energy of the rod will be given by
\begin{equation}
    \label{eq:energy_rod}
    E_{\rm rod}[\rm w, \Phi]:= E_{\rm el}[{\rm w}]+E_{\rm g}[{\rm w}]+E_{\rm sh}[\Phi].
\end{equation}
\subsection{Plateau problem: the spanning minimal surface}

Concerning the area contribution, the most recent approach is based on a new notion of contact between the film and the rod, named spanning. This has been introduced by Harrison and Pugh \cite{harrison2016existence} and later reformulated in \cite{DeLellis_2014} in a more Geometry Measure Theory context.
We refer also to \cite{Lussardi_2017, bevilacqua2019soap} for its application to the Kirchhoff-Plateau setting and to \cite{king2022plateau} for an interesting recent version of its use in the capillarity theory.
The spanning condition can be reformulated as follows
\begin{definition}
\label{cappio}
Let $H$ be a closed set in $\R^3$. We denote by $\mathcal{C}_H$ the set of all smooth embeddings $\gamma :S^1 \to \R^3 \setminus H$ which are not homotopic to a constant. Given a relatively closed subset $K$ in $\R^3 \setminus H$, we say that {\em $K$ spans $H$} if $\,\forall\, \gamma \in \mathcal{C}_H$, we have $K \cap \gamma (S^1) \neq \emptyset$.
\end{definition}
Then,  in \cite{DeLellis_2014}, for a fixed compact boundary, the following theorem has been proved.

\begin{theorem}\label{th:delellis}
Assume that there exists a surface $S \subset \R^3$ that spans $\R^3 \setminus \vect{p}[{\rm w}](\Omega)$ in the sense of Definition \ref{cappio} and such that $\mathcal{H}^2(S) < + \infty$. Then the following minimization problem has a solution, i.e. there exists
\begin{equation}
    \label{eq:ener_film}
    m_0 := {\rm min}\left\{\mathcal{H}^2(S): S \hbox{ spans } \vect{p}[{\rm w}](\Omega)\right\}.
\end{equation}
\end{theorem}
\begin{oss}
We remark that in Theorem \ref{th:delellis}, the position function $\vect{p}$ does not depend on $\Phi$ by the assumption that the boundary is a fixed compact set \cite{DeLellis_2014}. {\color{black}Moreover, the result has been carried out in the general framework of Almgren minimal sets \cite{almgren1968existence}.}
 \end{oss}

Hence, introducing the surface tension $\sigma>0$, the total energy of the system is given by\footnote{Denoting \textcolor{black}{the class of admissible smooth embeddings} with $F(\vect{p}[\rm w, \Phi](\Omega))$, as in Definition \ref{cappio}, we remark that 
$$
\inf_{{\rm w} \in U, S\in F(\vect{p}[{\rm w}, \textcolor{black}{\Phi}](\Omega))}\mathcal{E}[{\rm w}, S] = \inf_{{\rm w} \in U} \mathcal{E}_{\rm tot}[{\rm w}] = \inf_{{\rm w} \in U} \left(E_{\rm rod}[{\rm w}] + \inf_{S\in F(\vect{p}[{\rm w}, \textcolor{black}{\Phi}](\Omega))}E_{\rm f}[S]\right).
$$ }
\begin{equation}
    \label{eq:energy_tot}
   \mathcal{E}_{\rm tot}[{\rm w}, \Phi]:= E_{\rm el}[{\rm w}]+E_{\rm g}[{\rm w}]+E_{\rm sh}[\Phi] + 2 \sigma \,{\rm inf}\left\{\mathcal{H}^2(S): S \hbox{ spans } \vect{p}[{\rm w}, \textcolor{black}{\Phi}](\Omega)\right\}.
\end{equation}
\begin{oss}
	To prove the existence of a minimizer of \eqref{eq:energy_tot}, one needs to apply the Direct Method of the Calculus of Variations. However, such a proof is not an immediate consequence of Theorem 3.1 in \cite{Lussardi_2017} or Theorem 3.6 in \cite{bevilacqua2019soap}. The main difficulty is to show that the set of constraints $U$ \eqref{eq:set_constraint} is weakly closed in the right topology. In particular, to replicate Schuricht's result \cite{Schuricht_2002}, one has to introduce some conditions on the function $\Phi$ to take into account only some suitable deformations of the boundary. Since the aim of this paper is to provide some explicit examples of equilibria {\color{black} in a more regular setting} of the {\em modified} Kirchhoff-Plateau problem (penalizing shape modifications of the cross-section), we do not care about general results on the existence of minimizers of \eqref{eq:energy_tot}. Such a purpose will be a subject of a forthcoming work.
	\end{oss}

In the following, we consider a simplified version of the energy functional $\mathcal{E}_{\rm tot}$ in \eqref{eq:energy_tot}, due to the fact that the general problem is quite hard to study and to deal with explicit solutions is necessary to make some simplifications. Indeed, to provide a quantitative analysis of the minimizers and to make explicit calculations to understand the competition between the surface tension $\sigma$ and the elasticity of the rod, it is useful to set the problem in a more regular setting (for instance, the midline is going to be $\mathcal{C}^\infty$ and the soap film a regular parametrized surface). Moreover, we are forced to choose specific forms of the energy functionals in \eqref{eq:energy_tot} and to define a stronger spanning condition to know how the film spans the boundary of the rod in order to quantify boundary effects and normal forces applied  by the soap film to the rod.

\section{Plane midline energy functionals}
\label{sec:planar}

In this section, we aim to derive the Euler-Lagrange equations for a simplified version of the energy functional $\mathcal{E}_{\rm tot}$ introducing some additional hypotheses.

First of all, we require the midline $\vect{r}$ to be a \textcolor{black}{regular} plane closed unknotted curve lying in the plane $z = 0$. Hence, it can be parametrized in cartesian coordinates as $\vect{r}(s) = (x(s), y(s))$. The closure constraints can be just imposed on each component.

Second, \textcolor{black}{the curve is constrained to lie in the plane $z = 0$ so that the torsion of the midline vanishes. Third, we assume that there is no influence of the twist, i.e. $\omega=0$
}. 
In addition, we fix the bending parameters to be equal, namely $\kappa_1 = \kappa_2$, simplifying the elastic energy density. Indeed, we choose a specific energy density $f$ like $f(s,{\rm w})= \kappa(s)^2$ \textcolor{black}{with $\kappa^2=\kappa_1^2+\kappa_2^2$}. In cartesian representation, the elastic contribution becomes
$$
E_{\rm el}[\vect{r}]=\int_0^{L} (\kappa(s))^2 dl=\int_0^{L} \frac{(x'(s)y''(s)-x''(s)y'(s))^2}{((x'(s))^2+(y'(s))^2)^3} \sqrt{x'(s)^2+y'(s)^2} \,ds.
$$

Then, concerning the energy contribution of the film, we need to express in this setting the area of the surface $\mathcal{H}^2(S)$. Since the curve is fixed to be unknotted and planar, and no results are known about the regularity of the contact set \cite{DeLellis_2014}, we assume that the minimal surface attaches to the thick boundary in a curve in the plane $z = 0$. Hence, we introduce a function $a(s):[0,L]\to (0,+\infty)$ to measure the non-vanishing thickness of the cross-section $\mathcal{A}(s)$ and we require that it is small compared to the length of the curve. \textcolor{black}{Moreover, we require that $a(s)$ is differentiable to guarantee that the contact curve $\widetilde{\vect{r}}(s)$ is a regular curve on $[0,L]$. Finally, for the purpose of this paper and to characterize specific shapes of cross-sections, we also assume that $\mathcal{A}(s)$ is a convex set}. In such a setting, to make explicit calculations in evaluating the effects of surface tension $\sigma$ on different shapes of the cross-section $\mathcal{A}(s)$, we also require that the intersection between the film with the thick boundary is a plane smaller regular curve $\widetilde{\vect{r}}(s) = (u(s), v(s))$, called {\em contact curve}, whose expression, by a standard geometric argument, is the following one
\begin{align*}
&u(s) =x(s)-\frac{a(s)y'(s)}{\sqrt{(x'(s))^2+(y'(s))^2}},  &&v(s)=y(s)+\frac{a(s)x'(s)}{\sqrt{(x'(s))^2+(y'(s))^2}}.
\end{align*}
Hence, we can write the energy of the film as the area of the closed surface surrounded by the curve $\widetilde{\vect{r}}$ \cite[Section~I.6]{Elsgolc_2007}, namely
\begin{equation}
\label{Efilm}
\scalemath{0.8}{
\begin{aligned}
&E_{\rm f}[S]=2\sigma\mathcal{H}^2(S(u,v))=\sigma\int_0^{L} (u(s)v'(s)-v(s)u'(s))ds=\\ 
&=\sigma \int_0^{L} \left(x(s)-\frac{a(s)y'(s)}{\sqrt{(x'(s))^2+(y'(s))^2}}\right)\left(y'(s)+\frac{a'(s)x'(s)}{\sqrt{(x'(s))^2+(y'(s))^2}}+\frac{a(s)x''(s)}{\sqrt{(x'(s))^2+(y'(s))^2}}-\frac{a(s)x'(s)(x'(s)x''(s)+y'(s)y''(s))}{((x'(s))^2+(y'(s))^2)^{3/2}}\right)\,ds\\
&\quad-\sigma \int_0^{L} \left(y(s)+\frac{a(s)x'(s)}{\sqrt{(x'(s))^2+(y'(s))^2}}\right)\left(x'(s)-\frac{a'(s)y'(s)}{\sqrt{(x'(s))^2+(y'(s))^2}}-\frac{a(s)y''(s)}{\sqrt{(x'(s)^2+(y'(s))^2}}+\frac{a(s)y'(s)(x'(s)x''(s)+y'(s)y''(s))}{((x'(s))^2+(y'(s))^2)^{3/2}} \right)\,ds,
\end{aligned}
}
\end{equation}
where the surface tension $\sigma$ physically lies between $0$ and $1$.

Using a similar approach to the one employed in \cite{BBLM22}, to derive first-order necessary conditions for the energy functional, we move the fixed length constraint from the function space to the energy functional, introducing a Lagrange multiplier $\lambda \in \R$. The proof of the equivalence of the two approaches can be found in \cite{BBLM22}. So, to impose that the length of the curve is fixed, the additional functional term is the following
\begin{equation}
    \label{eq:inestensibilita}
    E_{\rm c}=\lambda \int_0^{L}\left(\sqrt{x'(s)^2+y'(s)^2}-1 \right)ds.
\end{equation}

By fixing the midline to remain in the plane $z = 0$, the weight contribution is irrelevant, hence in the derivation of the first-order necessary conditions we neglect it.

Now, the last term we introduce is the one which accounts for changes of shape of the cross-section $\mathcal{A}(s)$ to penalize big shape modifications. 
In \eqref{eq:E_sh_tot}, depending on the choice of $\Phi_{\vect{p}}$, we will change the type of sections we are going to study. To satisfy the unshearable constraint of the cross-section, we can consider either specific geometric shapes or affine transformations in the plane perpendicular to the tangent vector $\vect{t}$.

The specific examples of $\Phi_{\vect{p}}$ we are going to make in the following are:
\begin{enumerate}
    \item {\em elliptical cross-section with fixed area}: given a semi-axis of the ellipse \textcolor{black}{denoted with $a$} and the area of the ellipse \textcolor{black}{denoted with $\Pi$}, the form of $\Phi_{\vect{p}}$ is given by
    \begin{equation}
    \label{eq:sh_1}
         E_{\rm sh_{1}}[a] =\int_0^L \left(a(s)^2+\frac{\Pi^2}{\pi^2 a(s)^2}\right)\,dl;
    \end{equation}
    \item {\em elliptical section with \textcolor{black}{a dilation} of the horizontal semi-axis from the equilibrium configuration deduced in the previous case}: given the horizontal semi-axis $a_0$ of the ellipse in the equilibrium configuration of the previous case and $\Theta(s) >0 $ the \textcolor{black}{dilation} coefficient applied \textcolor{black}{at} each point of the cross-section $\mathcal{A}(s)$, the form of $\Phi_{\vect{p}}$ yields
     \begin{equation}
     \label{eq:sh_2}
         E_{\rm sh_{2}}[\Theta] = \int_0^{L} (\Theta(s)-1)^2 a_0^2\,dl;
     \end{equation}
    \item {\em oval cross-section with fixed area}: using as a  model for the oval a limacon of Pascal \cite{Lawrence_1972}, let $\Pi$ be the area of the oval and $a$ be one half of the sum between the shorter and the longer horizontal semi-axes of the oval, then $\Phi_{\vect{p}}$ reads
\begin{equation}
\label{eq:sh_3}
     E_{\rm sh_{3}}[a] = \int_0^{L} \left(\frac{2}{\pi}(\Pi-\pi a(s)^2)\right) \,dl.
\end{equation}
\end{enumerate}

\begin{oss}
Depending on the choice of $E_{\rm sh_{i}}$, we are considering a different shape of the cross-section $\mathcal{A}(s)$ appended on each point of the midline $\vect{r}$. Hence, in the reconstruction process of the rod to place it in $\R^3$, the expression of the position function $\vect{p}$ defined in \eqref{eq:defp} has to be modified according to the additional parameters $a$ or $\Theta$ \eqref{eq:sh_1} - \eqref{eq:sh_3}. 
\end{oss}



Since we consider a functional which is invariant under reparametrization, we adjust dimensional terms adding two positive constants $[\alpha]=$Nm$^2$ and $[\beta], [\lambda]=$N/m$^2$
. Hence, we end up with the functional $\mathcal{E}_{\rm tot_{i}}$ of the following form
\begin{empheq}{align}
\label{eq:functional}
\mathcal{E}_{\rm tot_{i}}[\vect{r},a,S,\lambda ]=\alpha E_{\rm el}[\vect
r]+E_{\rm f}[S]+E_{\rm c}+\beta E_{\rm sh_{i}}[a]=\int_0^{L} F(s,x,x',x'',y,y',y'',a,a',\lambda)\,\,ds.
\end{empheq} 

\begin{oss}
{\color{black} The energy functionals $E_{\rm sh_{i}}$ we are going to consider in Sections \ref{elsec}-\ref{sec:oval} are independent of $s$, namely the effects of the shape modification term is isotropic on each cross-section. This implies two facts. First of all, such a choice is compatible with the unshearability condition of the cross-section since $\mathcal{A}(s)$ behaves in the same way for all $s \in [0,L]$. Second, under such strict hypothesis, the existence of a minimizer for \eqref{eq:energy_tot} and in particular for \eqref{eq:functional}, in the planar setting, holds true in the weak topology. Indeed, the set of the constraints $U$ results to be weakly closed in $L^2$ and the Ciarlet-Ne\v cas condition holds true by the isotropic request of $E_{\rm sh_{i}}$ on each cross-section $\mathcal{A}(s)$ for all $s \in [0,L]$ \cite{Schuricht_2002}. Another issue is to remove such assumptions and looking for the weakest condition on admissible shape modifications of the cross-sections to get $U$ weakly closed in the right topology.}  
\end{oss}

{\color{black} Hence, we can} derive Euler-Lagrange equations for the energy functional $\mathcal{E}_{\rm tot_{i}}$ varying $E_{\rm sh_{i}}$ in the mentioned cases computing
\begin{empheq}[left=\empheqlbrace]{align}
\label{eq:prima_EL}
& \frac{d}{dx}F-\frac{d}{ds}\left(\frac{d}{dx'}F\right)+\frac{d^2}{ds^2}\left(\frac{d}{dx''}F\right)=0,\\
\label{eq:seconda_EL}
& \frac{d}{dy}F-\frac{d}{ds}\left(\frac{d}{dy'}F\right)+\frac{d^2}{ds^2}\left(\frac{d}{dy''}F\right)=0,\\
\label{eq:terza_eul_lagr}
& \frac{d}{da}F-\frac{d}{ds}\left(\frac{d}{da'}F\right)=0,\\
\label{eq:quarta_EL}
& \frac{d}{d\lambda}F=0.
\end{empheq}

In the following, we make a necessary assumption to characterize solutions of \eqref{eq:prima_EL} - \eqref{eq:quarta_EL}. 
Indeed, we aim to visualize at least a minimizer and to analyze quantitatively the interplay between elasticity and surface tension. Hence, since the system is quite complicated and there are many independent variables, it can be just faced using software for symbolic calculus, like Mathematica (Wolfram Inc., version 13). Moreover, to effectively solve \eqref{eq:prima_EL} - \eqref{eq:quarta_EL}, we choose an explicit form of the midline, substituting, once we have computed the derivatives, 
$$
(x(s), y(s))=\left(R \cos \left(\frac{2\pi s}{L}\right), R \sin \left(\frac{2\pi s}{L}\right)\right).
$$ 
Such a choice implies that we are always considering as midline a circumference of radius $R$ which has fixed length equal to $L = 2 \pi R$. This is an admissible solution among the critical ones for the Euler elastica in the plane \cite{langer1985curve}. \textcolor{black}{Assuming a constant surface tension, its action is homogeneous in every point of the surface and every direction, so the selection of a circular midline seems to be a reasonable choice in this simple isotropic context.}

\section{Case 1: Circular midline with elliptical section of fixed area}
\label{elsec}
Let $\Pi\in]0,+\infty[$ be the fixed area of the cross-section, $a(s)$ and $b(s)=\frac{\Pi}{\pi a(s)}$ be the semi-\textcolor{black}{axis} of the ellipse, where we choose $a(s)$ to be the semi-axis in the plane $z=0$ and $b(s) \perp a(s)$ for all $s \in [0,L]$.
The energy functional associated to the shape of the cross-section \eqref{eq:sh_1} can be rewritten as
\begin{equation*}
    E_{\rm sh_{1}}[a]=\int_0^{L} \left(a(s)^2+b(s)^2\right)\,dl= \int_0^{L} \left(a(s)^2+\frac{\Pi^2}{\pi^2 a(s)^2}\right) \sqrt{x'(s)^2+y'(s)^2} \,ds.
\end{equation*}

We report the derivation of the Euler-Lagrange equations in Appendix \ref{appendice1}. From \eqref{eq:terza_eul_lagr}, we can compute the value of $a$ as a root of the following quartic function
\begin{equation}\label{quarticaeq}
     \Gamma(a(s)):= \,\sigma \pi^2 a(s)^4+\beta \pi^2 Ra(s)^4-\pi^2 \sigma R a(s)^3-\beta \Pi^2 R=0.
\end{equation}
We notice immediately that the expression of $a(s)$ does not depend on $s$, implying that the action of the surface tension on the cross-section is the same on $[0,L]$.

By imposing the inextensibility constraint and since the midline is a circumference with fixed length $L$, \eqref{eq:sezione1} reduces to
\begin{empheq}[left=\empheqlbrace]{align}
& 2 \sigma R   -2 \sigma a + \beta a^2+\frac{\beta  \Pi^2}{\pi ^2 a^2} +\lambda -\frac{\alpha}{R^2}=0,\\
\label{eq:quartica2}
& \sigma \pi^2 a^4+\beta \pi^2 R a^4-\pi^2 \sigma R a^3-\beta \Pi^2 R=0.
\end{empheq}
\begin{figure}[htbp]
	\begin{subfigure}{.5\linewidth}
		\centering
		\includegraphics[width=0.94\textwidth]{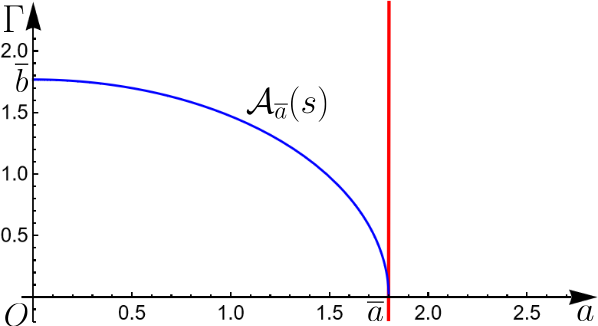}
		\caption{$\beta=1, \Pi=10, R=5$, $\sigma=0.1$}
		\label{quartica}
	\end{subfigure}
	\begin{subfigure}{.5\linewidth}
		\centering
		\includegraphics[width=\textwidth]{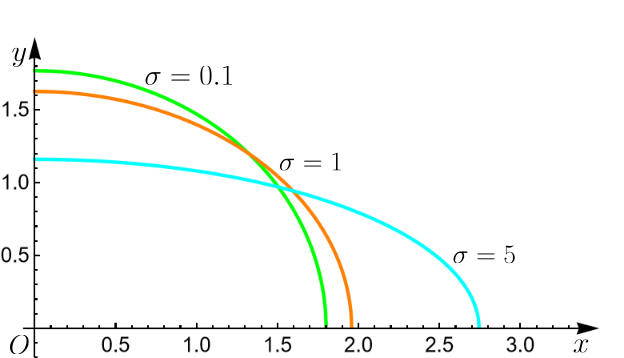}
		\caption{$\beta=1, \Pi=10, R=5$}
		\label{fig:vario_sigma}
	\end{subfigure}
	\caption{(a) The red line is the quartic $\Gamma(a)$ in \eqref{quarticaeq}. It has a unique positive real root $\overline{a}$. The blue line is the shape of the elliptical cross-section obtained at the stationary point $\overline{a}$. (b) Effect of the surface tension $\sigma$ on the elliptical shape of the cross-section. We choose $\sigma = 0.1, 1, 5$.}
	\label{fig:ellisse_caso1}
\end{figure}

Instead of solving explicitly \eqref{eq:quartica2} we study the behavior of solutions, computing 
the discriminant of \eqref{eq:quartica2} \cite{janson2010roots}, which reads
$$
\Delta = -256(\sigma \pi^2+\beta\pi^2 R)^3\Pi^6R^3-27\beta^2\pi^8\sigma^4\Pi^4R^6<0,
$$
implying that the equation has two distinct real roots and two complex conjugate ones. Next, we look for the real positive ones. Using the Descartes' rule \cite{cheriha2019descartes} for a fourth degree equation, we obtain that there exists a unique positive real root $\overline{a}$ for all possible choices of the parameters $(\beta, \sigma, \Pi, R)$, see Figure \ref{quartica}.


To verify that the obtained root $\overline{a}$ is a physical solution, we need to avoid interpenetration. Hence, we need to ensure $\overline{a} < R$. Computing the derivative of \eqref{eq:quartica2} with respect to $a$, we get
$$
\Gamma'(a)=4\sigma \pi^2 a^3+4\beta\pi^2 R a^3-3\pi^2\sigma R a^2=\pi^2a^2(4\sigma a+4\beta Ra-3R\sigma),
$$
which turns out to be positive if 
$$
a > \frac{3 \sigma R}{4\left(\beta R + \sigma\right)}=: a^{\star}.
$$
Since all constants are positive, we notice that $a^{\star} <R$. Assuming interpenetration of matter occurs, namely $\overline{a} \geq R$, it must be $\Gamma(R)\leq 0$, that results in $R \leq \sqrt{\frac{\Pi}{\pi}}$, which gives an explicit relation between the area of the section $\Pi$ and the radius of the midline $R$. Hence, the unique positive root $\overline{a}$ is admissible and physical if the radius of the circular midline satisfies the following inequality
\begin{equation}
\label{eq:non_comp}
    R > \sqrt{\frac{\Pi}{\pi}}.
\end{equation}

From \eqref{eq:non_comp}, we can also immediately establish if the ellipse is elongated along the horizontal axis or the vertical one. Indeed, since $\Pi=\pi \overline{a}\overline{b}$, $\Gamma(a)<0$ for all $a \in [0, \overline{a})$ and $\overline{a} < R$ by the non-interpenetration constraint, we necessarily have 
$$\Gamma\left(\sqrt{\frac{\Pi}{\pi}}\right)<0 \quad\Longrightarrow \quad\overline{a}>\left(\sqrt{\frac{\Pi}{\pi}}\right)\quad \Longrightarrow \quad
\overline{a}>\overline{b}.
$$



Finally, neglecting the effect of the surface tension $\sigma$, we get that the equilibrium cross-section $\overline{a}$ tends to be a circle, indeed $\overline{a} = \sqrt{\Pi/\pi} = \overline{b}$. Considering $\sigma \neq 0$ and increasing its value, see Figure \ref{fig:vario_sigma}, we notice that the elliptical shape of the cross-section elongates along the horizontal axis since it is pulled by the effect of the normal force generated on the boundary by the film \cite{roman2010elasto}.



\section{Case 2: Circular midline and elliptical section with \textcolor{black}{dilation} of the horizontal semi-axis}
\label{sec:homotetia}
\begin{figure}[htbp]
	\begin{subfigure}{.5\linewidth}
		\centering
		\includegraphics[width=0.7\textwidth]{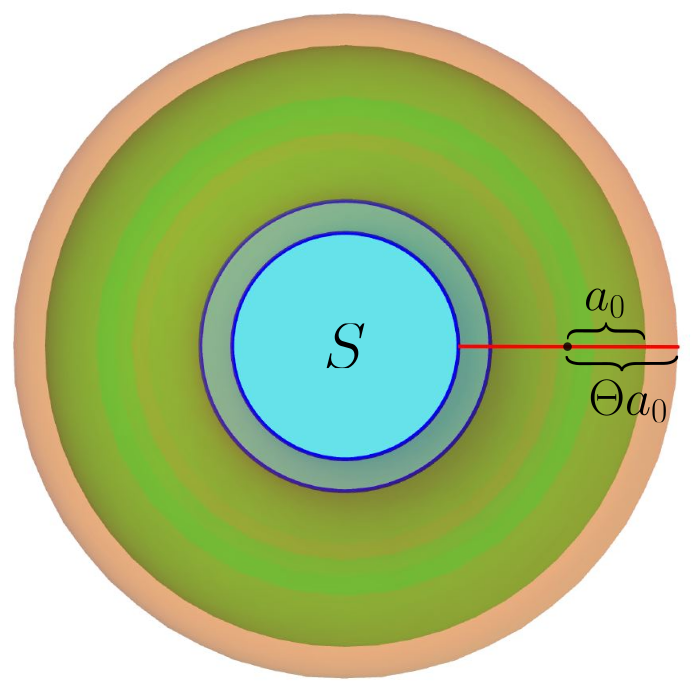}
		\label{torel1}
	\end{subfigure}
	\begin{subfigure}{.5\linewidth}
		\centering
		\includegraphics[width=0.77\textwidth]{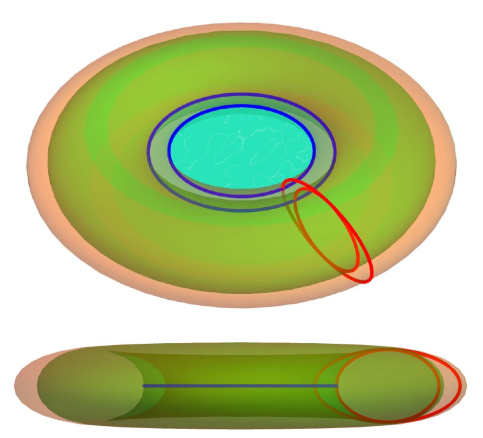}
		\label{torel2}
	\end{subfigure}
	\caption{(a) The system with circular midline and elliptical cross-section, viewed from above: in green we consider a torus with red elliptical cross-section (horizontal semi-axis $a_0$) (solution of Section \ref{elsec}); in orange, we display the torus with the \textcolor{black}{dilation} of the horizontal semi-axis, that is $\Theta a_0$, with $\Theta >0$.
(b) 3D visualization of the system and the elliptical section.}
	\label{torel}
\end{figure}

The starting point of this case is the equilibrium configuration obtained in Section \ref{elsec}: let us call the equilibrium shape of the cross-section $\overline{a}=:a_0$. In this setting, we consider a homothetic transformation of the horizontal semi-axis, namely $\Theta(s) a_0$, where $\Theta: [0,L]\to \mathbb{R}$ is a differentiable function that satisfies the closure condition, i.e. $\Theta(0)=\Theta(L)$ and $\Theta(s)>0$, see Figure \ref{torel}.

The energy functional \eqref{eq:sh_2} associated to the shape of the cross-section modified by a homothetic transformation can be rewritten as follows
\begin{align*}
   E_{\rm sh_{2}}[\Theta]
    =\int_0^{L} \norm{\begin{pmatrix} \Theta(s)-1 & 0\\ 0 & 0\end{pmatrix} \begin{pmatrix} a_0\\ b_0\end{pmatrix}}^2dl=\int_0^{L} (\Theta(s)-1)^2(a_0)^2 \sqrt{x'(s)^2 + y'(s)^2}\, ds.
\end{align*}


To derive the Euler-Lagrange equations, differently from Section \ref{elsec}, in \eqref{eq:terza_eul_lagr}, the derivation has been carried out with respect to the homothetic parameter $\Theta(s)$. We report all the calculations in Appendix \ref{appendice1}.
As before, the expression of the \textcolor{black}{dilation} parameter $\Theta$ can be found by \eqref{eq:terza_eul_lagr} and it does not depend on $s$. Since the circular midline has fixed length, \eqref{eq:sezione2} reduces to 
\begin{empheq}[left=\empheqlbrace]{align}
& \beta a_0^2  (\Theta-1)^2-2 a_0 \sigma \Theta  +\lambda -\frac{\alpha}{R^2}+2 R \sigma=0,\\
& 2 a_0 (a_0 \Theta (\beta  R+\sigma )-R (\beta a_0  +\sigma ))=0.
\end{empheq}
We obtain that the unique positive constant solution is given by
$$
\overline{\Theta}=\frac{R(\beta a_0+\sigma)}{a_0(\beta R+\sigma)}.
$$
Indeed, since $\overline{\Theta} > 1$, the associated homothety is a \textcolor{black}{dilation}. The strict inequality results by the non-interpenetration constraint $a_0 < R$ imposed in Section \ref{elsec}.
\begin{figure}[htbp]
	\begin{subfigure}{.5\linewidth}
		\centering
\includegraphics[width=\textwidth]{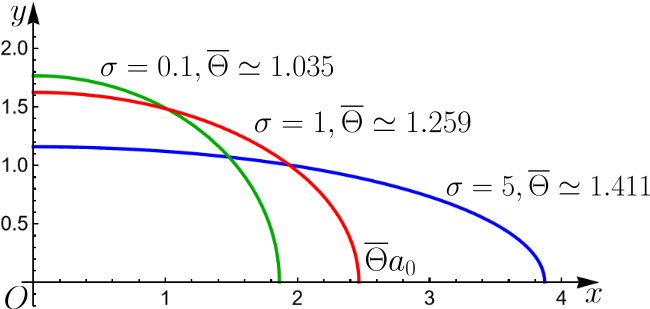}
		\caption{$\beta=1, R=5$}
\label{sigma_omega_bar}
	\end{subfigure}
	\begin{subfigure}{.5\linewidth}
		\centering
\includegraphics[width=0.98\textwidth]{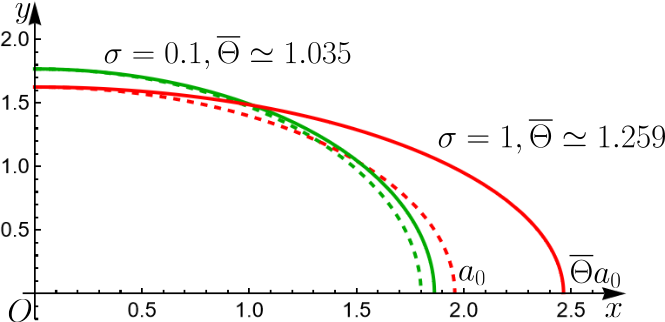}
		\caption{$\beta=1, R=5$}
\label{sigma_confronto}
	\end{subfigure}
	\caption{(a) Effect of the surface tension $\sigma$ on the elliptical shape of the cross-section modified via a \textcolor{black}{dilation} ($\Theta >1$). We choose $\sigma = 0.1, 1, 5$ and $a_0 \simeq 1.800, 1.958, 2.745$. (b) Effect of the surface tension $\sigma$ on the elliptical shape of the cross-section modified via a \textcolor{black}{dilation} $\Theta >1$ (continuum line) compared with the equilibrium solution $a_0$ obtained in Section \ref{elsec} (dashed line). We choose $\sigma = 0.1, 1$ and $a_0 \simeq 1.800, 1.958$.}
	\label{fig:omotetia}
\end{figure}
Finally, in Figure \ref{fig:omotetia}, we analyze the effect of the surface tension $\sigma$ on the shape of the cross-section. The presence of $\Theta$ results in a more elongated ellipse along the horizontal axis than the one obtained in Section \ref{elsec} by fixing the same value of the surface tension $\sigma$, see the dashed and the continuum curves in Figure \ref{sigma_confronto} .

\begin{oss}
We remark that we have considered an homothety only of the horizontal semi-axis $a$ for two reasons. On a hand, modifying the vertical one $b$ results in no modifications since the film spans the boundary along the horizontal axis of the cross-section. On the other hand, using two different coefficients to modify both $a$ and $b$ is equivalent to elongate/shorten $a$ while $b$ remains fixed.
\end{oss}


\section{Case 3: circular midline with oval section of fixed area}
\label{sec:oval}
\begin{figure}[ht]
	\begin{subfigure}{.5\linewidth}
		\centering
		\includegraphics[width=0.73\textwidth]{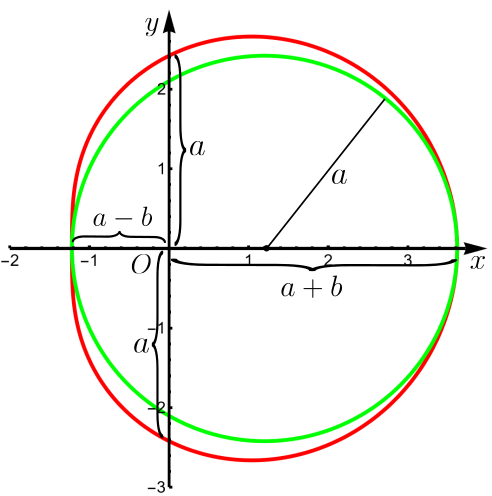}
		\caption{Red oval shape vs green circular one}
		\label{sezione_ovale}
	\end{subfigure}
	\begin{subfigure}{.5\linewidth}
		\centering
		\includegraphics[width=0.66\textwidth]{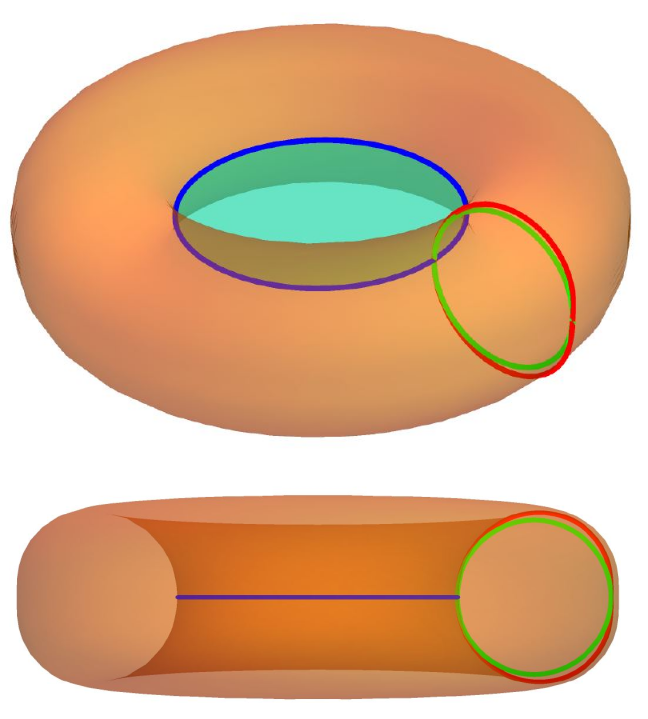}
		\caption{3D visualization of the oval shape}
		\label{3D_ovale}
	\end{subfigure}
	\caption{(a) Comparison between a circular shape (green) and the oval one (red) underlying the semi-axis $a+b, a, a, a-b$. (b) 3D representation of the rod with oval cross-section.}
	\label{fig:ovale_sezione}
\end{figure}

As in the first case, we fix the area $\Pi\in]0,+\infty[$ of the oval section and we parametrize it through the limacon of Pascal. In polar coordinates it has the form $r_\theta(s)=a(s)+b(s) \operatorname{cos}\theta$, with $a(s)>0$ and $b(s)\geq 0$. The vertical semi-axis is then $a(s)$, while the horizontal ones are $a(s)+b(s)$ and $a(s)-b(s)$, see Figure \ref{sezione_ovale}.\\ 

Moreover, to have a well-defined section without singularities, we first assume $a(s)>b(s)$, $\forall s\in[0,L]$, since the limacon has no double points if and only if $a(s)>b(s)$ \cite{Lawrence_1972}. Then, to avoid cusps, we require that $a(s)>2b(s)$, $\forall s\in[0,L]$.
Hence, the area of the oval is given by 
\begin{empheq}{align}
\label{area} \Pi=\int_0^{2\pi} \frac{1}{2}r_\theta(s)^2 d\theta=\int_0^{2\pi} \frac{1}{2}\left(a(s)^2+b(s)^2 \operatorname{cos} ^2\theta+2a(s)b (s)\operatorname{cos}\theta\right)d\theta=\pi a(s)^2+\frac{\pi}{2}b(s)^2,
\end{empheq}
which gives an explicit expression of $b$ as a function of $a$ and $\Pi$ as follows
\begin{equation}
\label{radicando}
 b(s)=\sqrt{\frac{2}{\pi}(\Pi-\pi a(s)^2)}.   
\end{equation}
The energy functional $E_{\rm sh_{3}}$ \eqref{eq:sh_3} can be then rewritten as
\begin{equation}
\label{eq:ESH_3}
    E_{\rm sh_{3}}[a]=
    \int_0^{L} \left(\frac{2}{\pi}(\Pi-\pi a(s)^2)\right) \sqrt{x'(s)^2+y'(s)^2} \,ds.
\end{equation}
We choose such an energy functional in order to penalize shapes too far from being a circle (\eqref{eq:ESH_3} is exactly $\int_0^L (b(s))^2\, dl$).

Now, we proceed similarly as done in Sections \ref{elsec} and \ref{sec:homotetia} deriving the Euler-Lagrange equations \eqref{eq:prima_EL} - \eqref{eq:quarta_EL} and substituting a circumference of radius $R$ for the midline $\vect{r}$ to understand the interaction between the elastic cross-section and the surface tension $\sigma$.
Differently from Section \ref{elsec}, however, here we require that the film spans the longer semi-axis of the oval, precisely $a+b$, see Figure \ref{3D_ovale}. We report the complete calculations for the derivation of \eqref{eq:prima_EL} - \eqref{eq:quarta_EL} in Appendix \ref{appendice1}.

From \eqref{eq:terza_eul_lagr}, we can compute the equilibrium solution $\overline{a}$ which is independent of $s$, as before.
Then, \eqref{eq:sezione6} reduces to
\begin{empheq}[left=\empheqlbrace]{align}
& 2 R^2 \left(\Pi \beta -\sqrt{2 \pi } \sigma  \sqrt{\Pi-\pi  a^2}\right)-\pi   \left(\alpha +R^2 (2 a (\beta  a+\sigma )-\lambda )\right)+2 \pi  R^3 \sigma =0,\\
\label{eq: terza ovale}
& \sqrt{2}\sigma  (\pi  a (2 a-R)-\Pi)+\sqrt{\Pi-\pi  a^2}\left(2\sqrt{\pi}\beta Ra+  \sigma \sqrt{\pi}(a + R)\right)=0.
\end{empheq}
In addition, such a solution $\overline{a}$ to be admissible must satisfy some constraints. First of all, the energy functional must be always well-defined, i.e., by \eqref{radicando},
 \begin{equation}
    \label{1.}
       \Pi-\pi a^2\geq 0, \quad \Longrightarrow \quad \overline{a}\leq\sqrt{\Pi/\pi}.
    \end{equation}
Second, to avoid cusps, we need to require $a > 2 b$, namely
\begin{equation}
    \label{2.} 
    a>2\sqrt{\frac{2}{\pi}(\Pi-\pi a^2)}, \quad \Longrightarrow\quad \overline{a}>\sqrt{\frac{8\Pi}{9\pi}}.
    \end{equation}
Finally, to avoid interpenetration of matter, 
$a+b=a+\sqrt{\frac{2}{\pi}(\Pi-\pi a^2)}<R$ results in defining a precise interval for the equilibrium solution $\overline{a}$:
\begin{equation}
    \label{3.}
     3\overline{a}^2-2R\overline{a}+R^2-2\frac{\Pi}{\pi}>0, \quad \Longrightarrow\quad \overline{a} < \frac{\pi R - \sqrt{2\pi(3 \Pi -  \pi R^2)}}{3 \pi} \,\cup \,\overline{a} > \frac{\pi R + \sqrt{2\pi(3 \Pi - \pi R^2)}}{3 \pi}.
    \end{equation} 
    

However, getting an explicit expression of $\overline{a}$ is quite hard. Hence, we look for a solution through a graphical method using as independent variables $a$ and $b$ which are independent of $s$. 
First of all, collecting the terms in \eqref{eq: terza ovale}, we obtain
\begin{equation}
 \label{hyperbole}
      2a^2-b^2+\left(1+\frac{2 \beta R}{\sigma}\right) ab-2Ra+Rb=0
\end{equation}
which is an hyperbola passing through the three points $O = (0,0), P_1 = (R,0)$ and $P_2 = (0,R)$. Second, from the constraint of the fixed area $\Pi$ of the oval \eqref{area}, rearranging the terms, we get an ellipse
\begin{equation}
\label{ellipse}
    2a^2+b^2=\frac{2\Pi}{\pi},
\end{equation}
with horizontal semi-axis equal to $\sqrt{\Pi/\pi}$. Then, we notice that $a$ must belong to the interval determined by \eqref{1.} and \eqref{2.}, and \eqref{3.} limits the range to $a+b<R$. Thus, we need to find a positive solution $(\overline{a}, \overline{b})$ which fulfils \eqref{hyperbole} - \eqref{ellipse} and lies in the region of the plane defined by \eqref{1.}, \eqref{2.} and \eqref{3.}.

\begin{figure}[htb]
\begin{subfigure}{.33\linewidth}
		\centering
\includegraphics[width=\textwidth]{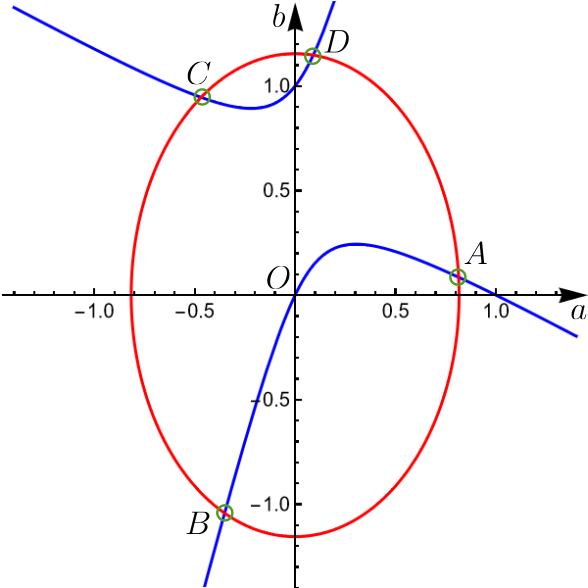}
\caption{$\beta=1, \sigma=0.9, R=1, \Pi=\frac{\pi}{3}$}
\label{fig:ell_hyp}
	\end{subfigure}
	\begin{subfigure}{.33\linewidth}
		\centering
\includegraphics[width=\textwidth]{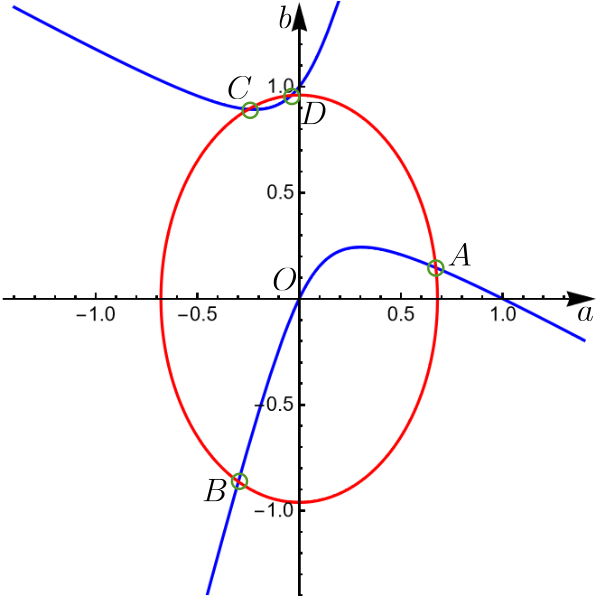}
	\caption{$\beta=1, \sigma=0.9, R=1, \Pi=\frac{29}{20}$}
  \label{fig:ell_hyp2}
	\end{subfigure}
 \begin{subfigure}{.33\linewidth}
		\centering
\includegraphics[width=\textwidth]{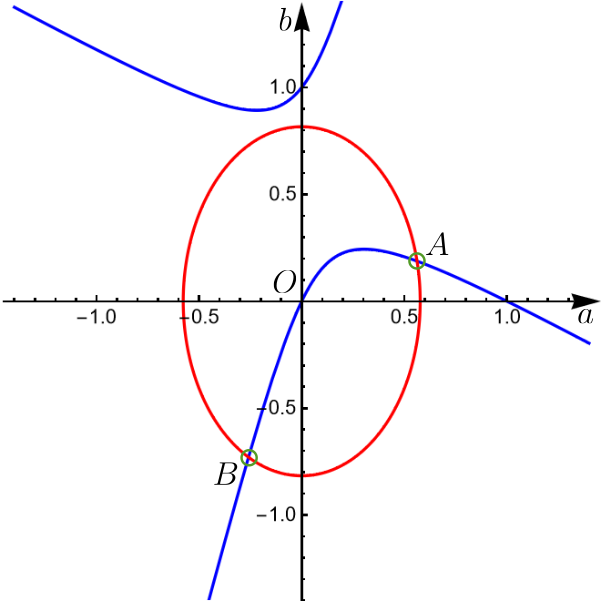}
\caption{$\beta=1, \sigma=0.9, R=1, \Pi=\frac{2\pi}{3}$}
  \label{fig:ell_hyp3}
	\end{subfigure}
	\caption{(a) Graphical representation of the ellipse (red) \eqref{ellipse} and the hyperbola (blue) \eqref{hyperbole} showing the intersection points fixing the parameters $R=1, \sigma=0.9, \beta=1$. 
 For $\Pi = \frac{\pi}{3}>\Xi_2 \simeq 1.3409$, we have four intersections. b) Graphical representation of \eqref{ellipse} and \eqref{hyperbole} with $\Pi = \frac{29}{20}>\Xi_2 \simeq 1.3409$, so again we obtain four intersections.  c) Graphical representation of \eqref{ellipse} and \eqref{hyperbole} with $\Pi = \frac{2\pi}{3}<\Xi_2 =1.3409$. Here, we have two intersections.
}
\label{fig:intersezioni}
\end{figure}
First of all, the intersections between the hyperbola \eqref{hyperbole} and the ellipse \eqref{ellipse} can be two or four. Fixing $R, \sigma$ and $\beta$, we have two points of intersection $A$ and $B$ if $\Pi\in (0,\Xi_2)$, four intersections $A, B, C$ and $D$ if $\Pi>\Xi_2$ and if $\Pi=\Xi_2$ the hyperbola and the ellipse are tangent (all the computations are reported in Appendix \ref{appendice2} together with the expression of $\Xi_2$, which is a root of the intersection polynomial).
In Figure \ref{fig:intersezioni}, we show different cases varying the value $\Pi$, i.e. the area of the cross-section, choosing suitable values for the parameters: $R= 1, \sigma= 0.9$ and $\beta=1$.

\begin{figure}[htb]
\begin{subfigure}{.5\linewidth}
		\centering
		\includegraphics[width=0.85\textwidth]{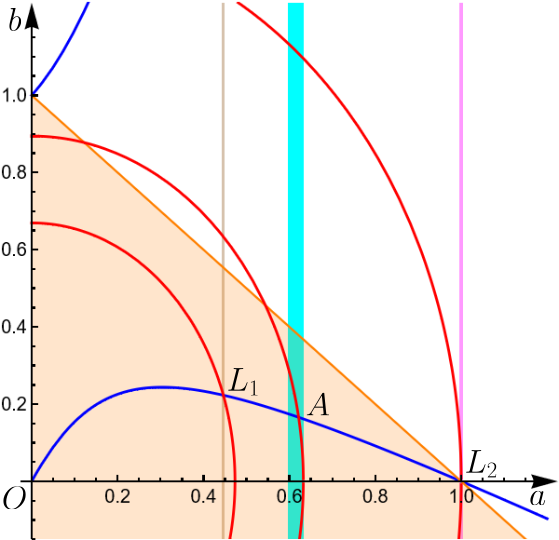}
		\caption{$\beta=1, \sigma=0.9,  R=1$}
		\label{fig:coniche}
	\end{subfigure}
	\begin{subfigure}{.5\linewidth}
		\centering
		\includegraphics[width=0.91\textwidth]{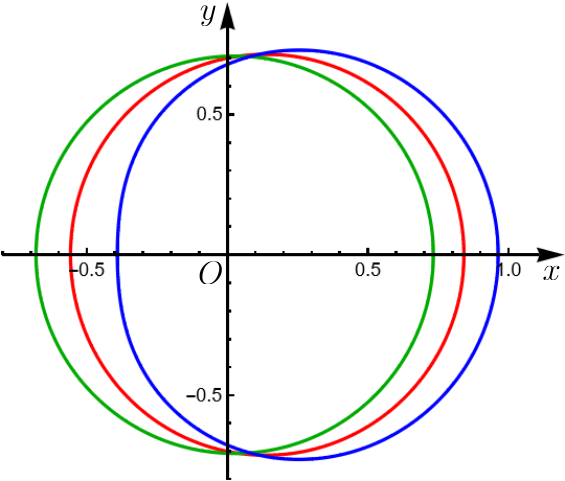}
		\caption{$\beta=1, \sigma=0.9, R=1$}
  \label{fig:ovalsection}
	\end{subfigure}
	\caption{(a) Graphical study of the real positive intersection point $A$ obtained fixing $R = 1$, $\sigma = 0.9$, $\beta = 1$ and $\Pi = \frac{2\pi}{5}$ and satisfying the non-interpenetration constraint \eqref{3.} represented by the orange half-plane. The solution belongs to the light blue interval obtained from \eqref{1.} and \eqref{2.}. Using the same parameters, we calculate $\tau\simeq 0.224$. Fixing $\Pi=\tau \pi R^2$,  we plot the first limit ellipse, that intersects the hyperbola and the brown line \eqref{2.} in the point $L_1$. Fixing $\Pi=\pi R^2=\pi$, we obtain the second limit ellipse, that intersects the hyperbola, the line $a+b=1$ and the pink line \eqref{1.} in the point $L_2=(1,0)$. b) Effect of the surface tension $\sigma$ on the oval cross-section. We fix $R=1, \beta=1, \Pi=\frac{\pi}{2}$ and we choose $\sigma=0.1$ (green line), $\sigma=1$ (red line) and $\sigma=10$ (blue line).}
 \label{limits}
 \end{figure}

We can immediately discard the intersection points with negative coordinates and we have to check if the remained ones satisfy \eqref{1.}, \eqref{2.} and \eqref{3.}. First of all, just by imposing \eqref{3.}, we conclude that there exists always (i.e. for each choice of the parameters $R, \sigma$ and $\beta$) a unique critical admissible intersection point $A:= (a(A), b(A))$. Then, we have to check if $A$ satisfies \eqref{1.} and \eqref{2.}.
In Figure \ref{fig:coniche}, we represent the plane \eqref{3.} in orange, the hyperbola \eqref{hyperbole} obtained fixing $R=1, \sigma=0.9, \beta=1$ and choosing $\Pi=\frac{2\pi}{5}$, we plot the ellipse \eqref{ellipse} that defines the point $A$. The light blue interval is obtained by \eqref{1.} and \eqref{2.} showing, by the graphical method, that the intersection point $A$, for such a choice of the physical parameters, is admissible. We should repeat a similar argument varying the parameters in physical ranges.

Moreover, we make some general considerations on the limit cases, namely when \eqref{1.}, \eqref{2.} and \eqref{3.} hold as an equality.
First of all, since $A$ belongs to the ellipse \eqref{ellipse} with horizontal semi-axis equal to $\sqrt{\Pi/\pi}$, then \eqref{1.} is immediately satisfied by the intersection point $A$.
Second, there are no solutions if the area of the cross-section $\Pi$ is higher with respect to the radius $R$ of the midline. Indeed, $A$ does not satisfy \eqref{3.} if the semi-axis of the ellipse is greater than $R$ (remember that the hyperbola \eqref{hyperbole} intersects the horizontal axis in $P_1=(R,0)$). Then, it must hold $\sqrt{\Pi/\pi}\leq R$. If the equality holds, then \eqref{hyperbole}, \eqref{ellipse} and \eqref{3.} meet in the point $L_2=(R,0)$: the corresponding solution has a circular section ($a=R$ and $b=0$), the film disappears and self-intersection occurs in $O$, see Figure \ref{fig:coniche}. 

Finally, from \eqref{2.}, we can deduce a lower bound for the range of parameters. Indeed, since $A$ must satisfy \eqref{2.}, we impose that the equality holds and we discuss the consequences. Hence, we obtain $a=\sqrt{8\Pi/(9\pi)}$ and by \eqref{ellipse},  $b=\sqrt{2\Pi/(9\pi)}$. Substituting the two values in \eqref{hyperbole}, we get a limit value for the area $\Pi$, i.e. 
$$
\frac{\Pi}{\pi R^2}=\frac{2}{\left(2+\frac{8\beta R}{9\sigma}\right)^2}=:\tau.
$$
Denoting with $L_1$ the intersection point between the ellipse obtained fixing $\Pi=\pi R^2 \tau$ and the hyperbola 
(see Figure \ref{fig:coniche}), such a solution corresponds to an oval cross-section with $a=2b$, which is not admissible. For lower values of $\Pi$, the section displays cusps since \eqref{2.} is no longer true.
Hence, choosing $\tau<\Pi/(\pi R^2)\leq 1$, then there exists a unique  physical critical point for the functional \eqref{eq:functional} with an oval cross-section with fixed area of the form \eqref{eq:ESH_3}.


Finally, in Figure \ref{fig:ovalsection}, we analyze the effect of the surface tension $\sigma$ on the shape of the cross-section. As expected, considering small values of $\sigma$, the section tends to be a circle (green line) since $\overline{b}\ll\overline{a}$. Increasing $\sigma$, we obtain oval sections with a more elongated horizontal semi-axis $a+b$, which is the one spanned by the soap film, while the other horizontal semi-axis $a-b$ is reduced.


\begin{oss}
We remark that in the complementary case, where we assume that the film attaches to the shortest horizontal semi-axis of the oval, i.e. $a-b$, the problem has no solutions. Indeed, there are no intersection points between the hyperbola and the ellipse in the admissible plane $a-b < R$, determined by the non-interpenetration condition, with both real positive coordinates $\overline{a}$ and $\overline{b}$.
\end{oss}


\section{Conclusions and final remarks}
\label{sec:conclusions}
The study of critical points for the Kirchhoff-Plateau problem is hard and not completely well-understood. Indeed, the presence of a highly constrained midline and a weak topological spanning condition makes difficult to derive the corresponding variations.\\
For all these reasons, we focused our attention on a simplified version of the problem, making some necessary assumptions to carry out explicit calculations. For instance, we required higher regularity for the curve, we considered only a planar situation and we specified the shape of the trace of the minimal surface on the boundary of the rod. We introduced an additional energy term to penalize shape modifications of the cross-sections appended to the midline and we studied the interaction between the surface tension of the soap film and the elastic contribution of the rod.\\
We derived Euler-Lagrange equations in three cases specifying the shape of the cross-section: elliptical cross-section with fixed area, elliptical section with \textcolor{black}{dilation} of the horizontal semi-axis from the equilibrium configuration and oval cross-section with fixed area.

We conclude that in all the studied cases and for a suitable range of the physical parameters, there exists a unique critical point, solution of the Euler-Lagrange equations, satisfying the imposed constraints in a simplified version of the problem. We also noticed that higher values of the surface tension $\sigma$ modify dramatically the shape of the cross-section.

Further developments will be devoted to analytically and numerically deriving a similar system of equations for critical points without making any assumptions or simplifications.
\textcolor{black}{For instance, a first direction to investigate is to overcome the limitation of the Kirchhoff rod. Indeed, it is a useful tool to make explicit calculations but it has limits from the physical point of view, we refer for instance to the hypotheses of unshearability and inextensibility of the midline. A possibility is to introduce a polyconvex energy functional to describe the rod as a {\em real} $3$D object. However, it is not so immediate due to the huge number of imposed constraints making hard to find the right topology to get compactness.
Another direction to follow is to employ numerical techniques to better analyze the effects of mechanical deformations. Such a question is still open again due to the highly constrained structure of the problem. For instance, we mention \cite{BBLM22}, where the authors provided some numerical tests, making some necessary simplifications and for a different version of the problem: they consider an elastic curve and not a rod. }

\section*{Acknowledgements}
The authors warmly thank Luca Lussardi and Alfredo Marzocchi for helpful suggestions and fruitful discussions.

GB is supported by the European Research Council (ERC), under the European Union's Horizon 2020 research and innovation program, through the project ERC VAREG - {\em Variational approach to the regularity of the free boundaries} (grant agreement No. 853404). GB is also supported by Gruppo Nazionale per l'Analisi Matematica, la Probabilità e le loro Applicazioni (GNAMPA) of Istituto Nazionale per l'Alta Matematica (INdAM).  CL is supported by Gruppo Nazionale per la Fisica Matematica (GNFM) of Istituto Nazionale per l'Alta Matematica (INdAM).


\begin{appendices}
\section{Expressions of Euler-Lagrange equations}
\label{appendice1}
In the following, we report the complete expressions of the Euler-Lagrange system in the three studied cases. We define $t:=\frac{2 \pi s}{L}$ to simplify the expressions.

In the first case, Section \ref{elsec}, where  we consider an ellipse as the cross-section and we substitute the circular midline, \eqref{eq:prima_EL} - \eqref{eq:quarta_EL} are given by
\begin{equation}
    \label{eq:sezione1}
   \left\{
   \scalemath{0.82}{
\begin{aligned}
&
-\frac{2\sigma \cos (t) a(s)a''(s)}{R}+\left(\left(-\frac{ \beta  \Pi^2  }{\pi ^2 a(s)^3}- \sigma  +a(s)\left(\beta  +\frac{\sigma}{R}\right)\right)\sin (t)-\frac{ \sigma  \cos (t) a'(s)}{R} \right)2a'(s)
\\ 
& +2a(s)\cos (t)\left(\frac{\beta  a(s) }{2}-\sigma \right)+\frac{\beta  \Pi^2 \cos (t)}{\pi ^2 a(s)^2}+\cos (t) \left(\lambda -\frac{\alpha }{R^2}+2 R \sigma \right)&=0,\\
& -\frac{2\sigma \sin (t) a(s)a''(s)}{R}+\left(\left(\frac{ \beta  \Pi^2  }{\pi ^2 a(s)^3}+ \sigma  -a(s)\left(\beta  +\frac{\sigma}{R}\right)\right)\cos (t)-\frac{ \sigma  \sin (t) a'(s)}{R} \right)2a'(s)+\\ 
& +2a(s)\sin (t)\left(\frac{\beta  a(s) }{2}-\sigma \right)+\frac{\beta  \Pi^2 \sin (t)}{\pi ^2 a(s)^2}+\sin (t) \left(\lambda -\frac{\alpha }{R^2}+2 R \sigma \right)&=0,\\
&  2a(s)\sigma -2 R \sigma-\frac{2 \beta  \Pi^2 R}{\pi ^2 a(s)^3}+2 a(s) \beta  R&=0,\\
& \frac{L}{2\pi R}-1 &= 0.
\end{aligned}
}
\right.
\end{equation}

In Section \ref{sec:homotetia}, we consider as cross-section an ellipse subjected to a \textcolor{black}{dilation} along the horizontal semi-axis and as midline the circumference. Hence \eqref{eq:prima_EL} - \eqref{eq:quarta_EL} simplify into

\begin{equation}
    \label{eq:sezione2}
    \left\{
    \scalemath{0.82}{
\begin{aligned}
& \frac{2\sigma a_0^2}{R}\cos(t) \Theta(s)\Theta''(s)+2a_0\Theta'(s) \left(\sin(t)\left(\beta a_0+\sigma-a_0\Theta(s)\left(\beta+\frac{\sigma}{R}\right)\right)+ \frac{\sigma a_0}{R}\cos (t)\Theta'(s)\right)+\\
& +2a_0\cos (t)\Theta(s)\left( \beta a_0+\sigma-\frac{1}{2}\beta a_0 \Theta(s)\right)-\cos (t) \left(\beta a_0^2+\lambda-\frac{\alpha}{R^2}+2 R \sigma\right)&=0,\\
& \frac{2\sigma a_0^2}{R}\sin(t) \Theta(s)\Theta''(s)+2a_0\Theta'(s) \left(\cos(t)\left(-\beta a_0-\sigma+a_0\Theta(s)\left(\beta+\frac{\sigma}{R}\right)\right)+ \frac{\sigma a_0}{R}\sin (t)\Theta'(s)\right)+\\
& +2a_0\sin (t)\Theta(s)\left( +\beta a_0+\sigma-\frac{1}{2}\beta a_0 \Theta(s)\right)-\sin (t) \left(\beta a_0^2+\lambda-\frac{\alpha}{R^2}+2 R \sigma\right)&=0,\\
& 2 a_0 (a_0 \Theta(s) (\beta  R+\sigma )-R (a_0 \beta +\sigma ))&=0,\\
& \frac{L}{2\pi R}-1 &=0.
\end{aligned}
}
\right.
\end{equation}

In the third case, Section \ref{sec:oval}, where  we consider an oval as the cross-section and a circular midline, \eqref{eq:prima_EL} - \eqref{eq:quarta_EL} become 
\begin{equation}
\label{eq:sezione6}
\left\{
\scalemath{0.82}{
\begin{aligned}
& 2 \sigma R a''(s)\cos (t) \sqrt{\Pi-\pi  a(s)^2} \left(\sqrt{2}  \left( \Pi^2   +2 \pi ^2  a(s)^4-3 \pi \Pi    a(s)^2\right)-\sqrt{\pi}a(s)\left(\Pi -\pi   a(s)^2\right)^{\frac{3}{2}}\right)+\\
& -2 \pi ^{3/2} R a(s)^2 a'(s)\left(\Pi-\pi  a(s)^2\right)\left( \sin (t) \left(a(s) (\sigma +2  \beta R  )+\sigma R \right) + \sigma   \cos (t) a'(s)  \right)+\\
&  +2\sqrt{2} \sigma R \sqrt{\Pi-\pi  a(s)^2}a'(s) \left(-   \sin (t) (\Pi^2+2  \pi ^2   a(s)^4-\pi ^2 R     a(s)^3)+2 \pi ^2  \cos (t) a(s)^3  a'(s)  \right)+\\
& +2\sqrt{2} \sigma R \sqrt{\Pi-\pi  a(s)^2}a'(s) \left(\pi\Pi a(s)( \sin (t)   (3a(s)-R)-3     \cos (t) a'(s))\right)+\\
&+2\Pi R \sqrt{\pi}a'(s)\left(\Pi-\pi  a(s)^2\right)\left(\sigma (-  \cos (t) a'(s) +R   \sin (t))+ \sin (t) a(s)(\sigma+2 \beta R)\right)+\\
& +\sqrt{\pi} \cos (t) \left(\Pi-\pi  a(s)^2\right)\left(\pi a(s)^2\left(-2  R^2  a(s) \left(\beta   a(s)+\sigma   \right)- \alpha   + \lambda R^2    +2 \sigma R^3  \right)\right)+\\
 & +\sqrt{\pi} \cos (t) \left(\Pi-\pi  a(s)^2\right)\left(+2 \Pi R^2 a(s)(2 \beta   a(s)+ \sigma )+\Pi(\alpha- \lambda R^2   -2 \sigma R^3 )\right)+\\
& +2\sqrt{2} \sigma R^2 \sqrt{\Pi-\pi  a(s)^2}\cos (t)\left(\pi  a(s)^2(\pi  a(s)^2-2  \Pi  )+ \Pi^2 \right)-\frac{2 \beta\Pi^2 R^2   \left(\Pi-\pi  a(s)^2\right) \cos (t)}{\sqrt{\pi }} & =0,\\
& 2\sigma R  a''(s)\sin (t) \sqrt{\Pi-\pi  a(s)^2} \left(\sqrt{2}  ( \Pi^2   +2 \pi ^2  a(s)^4-3 \pi \Pi    a(s)^2)-\sqrt{\pi}a(s)\left(\Pi-\pi   a(s)^2 \right)^{\frac{3}{2}}\right)+\\
& -2 \pi ^{3/2} R a(s)^2 a'(s)\left(\Pi-\pi  a(s)^2\right)\left( -\cos (t) \left(a(s) (2 \beta R    + \sigma) +\sigma R \right)+ \sigma   \sin (t) a'(s)   \right)+\\
&  +2\sqrt{2} \sigma R \sqrt{\Pi-\pi  a(s)^2}a'(s) \left(-   \cos (t) (-\Pi^2-2  \pi ^2   a(s)^4+\pi ^2 R     a(s)^3)+2 \pi ^2  \sin (t) a(s)^3 a'(s)  \right)+\\
& +2\sqrt{2} \sigma R \sqrt{\Pi-\pi  a(s)^2}a'(s) \left(\Pi\pi a(s)( \cos (t)   (-3a(s)+R)-3     \sin (t) a'(s))\right)+\\
&-2\Pi R\sqrt{\pi}a'(s)\left(\Pi-\pi  a(s)^2\right)\left(\sigma  (\sin (t) a'(s) +R   \cos (t))+ \cos (t)  a(s)(2R \beta  +\sigma)\right)+\\
& +\sqrt{\pi} \sin (t) \left(\Pi-\pi  a(s)^2\right)\left(\pi a(s)^2\left(-2  R^2  a(s) \left(\beta   a(s)+\sigma   \right)- \alpha   + \lambda  R^2   +2 \sigma R^3  \right)\right)+\\
 & +\sqrt{\pi} \sin (t) \left(\Pi-\pi  a(s)^2\right)\left(+2 \Pi R^2 a(s)(2 \beta   a(s)+ \sigma )+\Pi(\alpha- \lambda   R^2  -2 \sigma R^3 )\right)+\\
& +2\sqrt{2} \sigma R^2 \sqrt{\Pi-\pi  a(s)^2}\sin (t)\left(\pi  a(s)^2(\pi  a(s)^2-2  \Pi  )+ \Pi^2 \right)-\frac{2 \beta \Pi^2 R^2   \left(\Pi-\pi  a(s)^2\right) \sin (t)}{\sqrt{\pi }} & =0,\\
& \frac{\sqrt{2}\sigma  (\pi  a(s) (2 a(s)-R)-\Pi)}{\sqrt{\Pi-\pi  a(s)^2}}+2\sqrt{\pi}\beta Ra(s)+  \sigma \sqrt{\pi}(a(s) + R)&=0,\\
& \frac{L}{2\pi R}-1 &=0.
\end{aligned}
}
\right.
\end{equation}
\section{Intersections between \texorpdfstring{\eqref{hyperbole}}{(6.9)} and \texorpdfstring{\eqref{ellipse}}{(6.10)}}
\label{appendice2}
To determine the number of intersections between the hyperbola \eqref{hyperbole} and the ellipse \eqref{ellipse}, we study the discriminant $\Delta$ \cite{janson2010roots} of the quartic equation obtained from the intersection system of the two conics, i.e.
$$
\begin{aligned}
(18  \sigma ^2 +8  \beta ^2 R^2 +8  \beta  R \sigma) \pi ^2 a^4
&+(-12 \pi ^2 R \sigma ^2 +8 \pi ^2 \beta  R^2 \sigma ) a^3
+(6 \pi ^2 R^2 \sigma ^2-8 \pi  \Pi \beta ^2 R^2  -18 \pi  \Pi \sigma ^2-8 \pi  \Pi \beta  R \sigma ) a^2\\
&+(-8   \Pi \beta  R^2 \sigma  +4  \Pi R \sigma ^2) \pi a
+(4 \Pi^2-2 \pi  \Pi R^2) \sigma ^2=0
\end{aligned}
$$
which can be solved with respect to $a$ and any solution depends on the physical parameters $R, \sigma, \beta, \Pi$. Since all the coefficients are real, if $\Delta >0$ the obtained four roots or are all real or all complex conjugate. Solving $\Delta = 0$ with respect to the area $\Pi$, we end up with $6$ roots $\Xi_i$ with $i = 1,\dots,6$  which are functions of $R, \beta, \sigma$ and they are given by
$$
\begin{aligned}
    &\Xi_1=0, &&\Xi_2>0, &&&\Xi_3= \Xi_4 =\frac{\pi  R^2 \sigma  (2 \beta  R+3 \sigma )}{(2 \beta  R+\sigma )^2}>0, &&&&\Xi_5, \Xi_6 \,\text{ (conjugate roots)},
\end{aligned}
$$
where
\begin{equation}
\scalemath{0.72}{
\begin{aligned}
& \Xi_2=\frac{\sqrt[3]{\delta+\sqrt{4 \left(-81 \pi ^2 R^4 \left(4 R^2 \beta ^2+4 R \sigma  \beta +9 \sigma ^2\right)^4 \sigma ^4-81 \pi ^2 R^4 \left(4 R^2 \beta ^2+4 R \sigma  \beta +9 \sigma ^2\right)^3 \left(4 R^2 \beta ^2-4 R \sigma  \beta -9 \sigma ^2\right) \sigma ^4\right)^3+\delta^2}}}{3 \sqrt[3]{2} \left(4 R^2 \beta ^2+4 R \sigma  \beta +9 \sigma ^2\right)^3} +\frac{3 \pi  R^2 \sigma ^2}{4 R^2 \beta ^2+4 R \sigma  \beta +9 \sigma ^2}+\\
& -\frac{\sqrt[3]{2} \left(-81 \pi ^2 R^4 \left(4 R^2 \beta ^2+4 R \sigma  \beta +9 \sigma ^2\right)^4 \sigma ^4-81 \pi ^2 R^4 \left(4 R^2 \beta ^2+4 R \sigma  \beta +9 \sigma ^2\right)^3 \left(4 R^2 \beta ^2-4 R \sigma  \beta -9 \sigma ^2\right) \sigma ^4\right)}{3 \left(4 R^2 \beta ^2+4 R \sigma  \beta +9 \sigma ^2\right)^3 \sqrt[3]{\delta+\sqrt{4 \left(-81 \pi ^2 R^4 \left(4 R^2 \beta ^2+4 R \sigma  \beta +9 \sigma ^2\right)^4 \sigma ^4-81 \pi ^2 R^4 \left(4 R^2 \beta ^2+4 R \sigma  \beta +9 \sigma ^2\right)^3 \left(4 R^2 \beta ^2-4 R \sigma  \beta -9 \sigma ^2\right) \sigma ^4\right)^3+\delta^2}}},
\end{aligned}
}
\end{equation}
and
\begin{equation}
\scalemath{0.75}{
\begin{aligned}
 \delta & =17915904 \pi ^3 \beta ^{12} \sigma ^6 R^{18}+89579520 \pi ^3 \beta ^{11} \sigma ^7 R^{17}+380712960 \pi ^3 \beta ^{10} \sigma ^8 R^{16}+985374720 \pi ^3 \beta ^9 \sigma ^9 R^{15}+2205895680 \pi ^3 \beta ^8 \sigma ^{10} R^{14}+\\
& +3545109504 \pi ^3 \beta ^7 \sigma ^{11} R^{13}+4963265280 \pi ^3 \beta ^6 \sigma ^{12} R^{12}+4988459520 \pi ^3 \beta ^5 \sigma ^{13} R^{11}+4336558560 \pi ^3 \beta ^4 \sigma ^{14} R^{10}+\\
& +2295825120 \pi ^3 \beta ^3 \sigma ^{15} R^9+1033121304 \pi ^3 \beta ^2 \sigma ^{16} R^8.
\end{aligned}
}
\end{equation}
Hence, for $\Delta >0$, we have an admissible solution if $\Pi \geq \Xi_2$, resulting in four intersection points, while if $\Delta < 0$, i.e. $\Pi \in (0, \Xi_2)$, the hyperbola and the ellipse have just two admissible intersections. 
\end{appendices}

\printbibliography

\end{document}